\documentclass[11pt,a4paper]{article}
\usepackage{bbm}
\usepackage{mathrsfs}
\usepackage{amsfonts}
\usepackage{}

\usepackage{graphicx}
\usepackage{amsfonts,amsmath, amssymb}
\usepackage{color}

\newtheorem{theo}{\textbf{\ \ \quad Theorem}}[section]

\newtheorem{remark}{\textbf{\ \ \quad Remark}}[section]

\newcommand{\lbl}[1]{\label{#1}}

\newcommand{\be}{\begin{equation}}
\newcommand{\ee}{\end{equation}}
\newcommand\bes{\begin{eqnarray}}
\newcommand\ees{\end{eqnarray}}
\newcommand{\bess}{\begin{eqnarray*}}
\newcommand{\eess}{\end{eqnarray*}}

\newcommand{\nm}{\nonumber}

{\setlength\arraycolsep{2pt}}

\title{Blowup solutions of Grushin's operator}

\author{Guangying Lv$^a$, Jinlong Wei$^b$ and Longjie Xie$^c$\\
\\
\ \\
   {\small \it $^a$Institute of Applied Mathematics, Henan University}\\
  {\small \it Kaifeng, Henan 475001, China}\\
  {\small \tt gylvmaths@henu.edu.cn}\\
  {\small \it $^b$ School of Statistics and Mathematics, Zhongnan University of}\\
  {\small \it
 Economics and Law, Wuhan, Hubei 430073, China}\\
   {\small \tt  weijinlong.hust@gmail.com }\\
  {\small \it $^c$ School of Mathematics and Statistics, Jiangsu Normal University}\\
  {\small \it
Xuzhou,
Jiangsu 221000, P.R.China}\\
   {\small \tt  xlj.98@whu.edu.cn}
}

\begin{document}
\maketitle

\medskip

\begin{abstract} In this note, we consider the blowup phenomenon of Grushin's operator.
By using the knowledge of probability, we first get expression of heat kernel of Grushin's operator.
Then by using the properties of heat kernel and suitable auxiliary function, we get that the solutions
will blow up in finite time.

{\bf Keywords}: Grushin's operator; Heat kernel; Blowup.

\textbf{AMS subject classifications} (2010): 35K20, 60H15, 60H40.

\end{abstract}

\baselineskip=15pt

\section{Introduction}
\setcounter{equation}{0}

The finite time blowup phenomenon has been studied by many authors, see the book \cite{Hubook2018}.
There are two cases to study this problem. One is bounded domain and the other is whole space.
In this paper, we only consider the problem in the whole space.
For the whole space, the following "Fujita Phenomenon" has been attraction
in the literature. Consider the following Cauchy problem
   \bes\left\{\begin{array}{lll}
u_t=\Delta u+u^p,\ \ \ &x\in\mathbb{R}^d, \ \ t>0,\ \ p>0,\\
u(0,x)=u_0(x), \ \ \ &x\in\mathbb{R}^d.
   \end{array}\right.\lbl{1.1}\ees
It has been proved that:
 \begin{quote}
(i) if $0<p<1$, then every nonnegative solution is global, but not necessarily
unique;

(ii) if $1<p\leq1+\frac{2}{d}$, then any nontrivial, nonnegative solution
blows up in finite time;

(iii) if $p>1+\frac{2}{d}$, then $u_0\in\mathcal{U}$ implies that
$u(t,x,u_0)$ exists globally;

(iv) if $p>1+\frac{2}{d}$, then $u_0\in\mathcal{U_1}$ implies that
$u(t,x,u_0)$ blows up in finite time,
 \end{quote}
where $\mathcal{U}$ and $\mathcal{U_1}$ are defined as follows
   \bess
\mathcal{U}&=&\left\{v(x)|v(x)\in BC(\mathbb{R}^d,\mathbb{R}_+),
v(x)\leq \delta e^{-k|x|^2},\ k>0,\delta=\delta(k)>0\right\},\\
\mathcal{U_1}&=&\left\{v(x)|v(x)\in BC(\mathbb{R}^d,\mathbb{R}_+),
v(x)\geq c e^{-k|x|^2},\ k>0,c\gg1\right\}.
  \eess
Here $BC=\{$ bounded and uniformly continuous functions $\}$, see
Fujita \cite{F1966,F1970} and Hayakawa \cite{kH1973}. The proof of case (i)-(iii) relies on
the properties of heat kernel and suitable auxiliary function. Comparison principle is the
main tool to prove case (iv). In this note, we consider the degenerate parabolic operator--Grushin's operator.
We will consider the first three cases.

There are a lot of known results about the blowup phenomenon of parabolic equations.
Blowup phenomenon of quasilinear parabolic equations with Robin boundary condition
was considered by Enache \cite{E2011}, also see \cite{E1999,LW2004}. Then the
blowup phenomena of degenerate parabolic and nonlocal diffusion equations
were considered by \cite{KT2006,LW2012,LC2017,MS2001,PSZ2016}. Seki \cite{S2018}
obtained the type II blowup mechanisms. Zhang-Wang \cite{ZW2018} considered the blowup phenomenon of
3-D primitive equations of oceanic and atmospheric dynamics.

In this note, we consider a special degenerate parabolic operator--Grushin's operator. Fortunately, we can obtain expression of Grushin's operator. In next section, some
preliminaries are given and the main results will be proved in section 3.
Throughout this paper, we write $C$ as a general positive constant and $C_i$, $i=1,2,\cdots$ as
a concrete positive constant.

\section{Main results}
\setcounter{equation}{0}
Consider the Grushin's operator
   \bess
\mathcal{L}=\frac{1}{2}(\partial_{x_1}^2+x^2\partial_{x_2}^2),
   \eess
which is the generator of the diffusion process $(X^1_t,X^2_t)$, where $(X^1_t,X^2_t)$ satisfies
  \bess\left\{\begin{array}{llll}
dX^1_t=dW^1_t,\\
dX^2_t=X_tdW^2_t,\\
X^1_0=\mu_1,\ \ \ X^2_0=\mu_2.
  \end{array}\right.\eess
Here $W^i_t$ is a  standard Brownian motion, $i=1,2$. It is easy to see that
the process $(X^1_t,X^2_t)$ is a Gaussian stochastic process. Direct calculations
show that
   \bess
\mathbb{E}\left(\begin{array}{cccc}  X^1_t\\[-.5mm]  X^2_t
\end{array}\right)=\left(\begin{array}{cccc}  \mu_1\\[-.5mm]  \mu_2
\end{array}\right),\ \ \
Cov(X^1_t,X^2_t)=\left(\begin{array}{cccc}
t &\ \mu_1t & \\[-.5mm]
\mu_1t &\ \mu_1^2t+\frac{1}{2}t^2 &
\end{array}\right) .
   \eess
Therefore, we get the heat kernel of the
operator $\mathcal{L}$ is
  \bess
K(t,x_1,\mu_1,x_2,\mu_2)=\frac{1}{2\pi t^{3/2}}\exp\left\{-\frac{(x_1-\mu_1)^2}{t}-\frac{[\mu_1(x_1-\mu_1)-x_2+\mu_2]^2}{t^2}\right\},
   \eess
which yields that
  \bess
\nabla_{x_1}K(t,x_1,x_2)=-\frac{2x_1}{t}K(t,x_1,x_2),\ \
\nabla_{x_2}K(t,x_1,x_2)=-\frac{y}{t^2}K(t,x_1,x_2).
  \eess
It is easy to see that for classical heat kernel, we have
 $x\sim\sqrt{t}$. But in our case, different axis has different
scaling, that is,
   \bess
x_1\sim\sqrt{t},\ \ \ x_2\sim t.
  \eess
Now, we consider the following degenerate parabolic equation
   \bes\left\{\begin{array}{lll}
u_t=\mathcal{L} u+u^p,\ \ \ &x\in\mathbb{R}^2, \ \ t>0,\ \ p>0,\\
u(0,x)=u_0(x)\geq0, \ \ \ &x\in\mathbb{R}^2.
   \end{array}\right.\lbl{2.1}\ees
The main results is as followings.
  \begin{theo}\lbl{t2.1} Assume that $u_0$ is a bounded continuous non-negative function.

  (i) If $0<p<1$, then the solution of (\ref{2.1}) exists globally.

(ii) If $1<p\leq1+\frac{2}{3}$, then all nontrivial solutions of
(\ref{2.1}) blow up in finite time. That is to say, there exists a positive $T^*>0$ such that
   \bess
u(t,x)=\infty, \ \ t>T^*.
  \eess

  (iii) If $p>1$, then the solution of
(\ref{2.1}) blows up in finite time provided the initial datum satisfies
  \bess
\inf_{x\in\mathbb{R}^2}u_0(x)\geq\mu>0,
  \eess
where $\mu$ is a constant.
  \end{theo}
\begin{remark}\lbl{r2.1}
Comparing with the classical parabolic equation, that is to say, comparing
\cite[Theorem 5.5]{Hubook2018} with the above theorem \ref{t2.1}, we find
the value of $p$ is different. More precisely, it follows \cite[Theorem 5.5]{Hubook2018} that
when $1<p\leq \frac{2}{d}$ ($d$ is the dimension of space), the solutions of (\ref{2.1}) with $\mathcal{L}$ replaced by $\Delta$ will
blow up in finite time under the condition that the initial data $u_0\geq0(\not\equiv0)$ is bounded continuous function.
However, in the case of (\ref{2.1}), the index is $1<p\leq1+\frac{2}{3}$.

The assumption of (iii) is too strict, one can weaken the assumption.
  \end{remark}


\section{Proof of main results}
\setcounter{equation}{0}

{\bf Proof of Theorem \ref{t2.1}}  The solution of (\ref{2.1}) can be expressed as
  \bes
u(t,x)=K\ast u_0(x)+\int_0^t\int_{\mathbb{R}^2}K(t-s,x-y)u^p(s,y)dy.
   \lbl{3.1}\ees
Due to the positivity of heat kernel, it is easy to see that if the initial data is non-negative, then
the solution will keep positive.
The equality (\ref{3.1}) yields that for any $T>0$ and  $t\in[0,T]$,
  \bess
u(t,x)\leq \sup_{x\in\mathbb{R}^2}|u_0(x)|+[\sup_{x\in\mathbb{R}^2}|u(t,x)|]^p,
   \eess
where we used the properties of heat kernel. Hence we have for any $T>0$ and  $t\in[0,T]$,
  \bess
\sup_{x,y\in\mathbb{R}}|u(t,x,y)|\leq C(T),
   \eess
which implies the result of (i). Denote $u(t,x)=I_1(x,t)+I_2(x,t)$ and
   \bess
I_1(t,x)=\int_{\mathbb{R}^2}K(t,x-y)u_0(y)dy,\ \
I_2(t,x)=\int_0^t\int_{\mathbb{R}^2}K(t-s,x-y)u^p(s,y)dyds.
   \eess
We may assume without loss of generality that $u_0(y)\geq C_1>0$ for $|y|<1$ by the
assumption. A direct computation shows that
   \bes
I_1(t,x)&\geq&\frac{C_1}{t^{3/2}}\int_{B_1(0)}\exp\left(-\frac{2x_1^2+2y_1^2}{t}-\frac{2x_2^2+2y_2^2}{t^2}\right)dy\nm\\
&\geq& \frac{C_1}{t^{3/2}}\exp\left(-\frac{2x_1^2}{t}-\frac{2x_2^2}{t^2}\right)\int_{|y|\leq\frac{1}{\sqrt{t}}}
\exp\left(-2y_1^2-\frac{2y_2^2}{t}\right)dy\nm\\
&\geq& \frac{C_1}{t^{3/2}}\exp\left(-\frac{2x_1^2}{t}-\frac{2x_2^2}{t^2}\right)
  \lbl{3.2} \ees
for $t>1$ and $C_1>0$.

It is easy to see that
   \bess
I_2(t,x)&\geq& \int_0^t\left(\int_{\mathbb{R}^2}K(t-s,x-y)u(s,y)dy\right)^pds.
   \eess
Let
   \bess
G(t)=\int_{\mathbb{R}^2}K(t,x) u(t,x)dx.
   \eess
Then for $t>1$,
   \bes
G(t)&=&\int_{\mathbb{R}^2}I_1(t,x)K(t,x)dx+\int_{\mathbb{R}^2}I_2(t,x)K(t,x)dx\nm\\
&\geq&\frac{C_2}{t^{3/2}}+\int_0^t\int_{\mathbb{R}^2}K(t,x)\left(\int_{\mathbb{R}^2}K(t-s,x-y)u(s,y)dy\right)^pdxds\nm\\
&\geq&\frac{C_2}{t^{3/2}}+\int_0^t\left[\int_{\mathbb{R}^2}\left(\int_{\mathbb{R}^2}K(t,x)K(t-s,x-y)dx\right) u(s,y)dy\right]^pds.
   \lbl{3.3}\ees
It is clear that
   \bess
&&\int_{\mathbb{R}^2}K(t,x)K(t-s,x-y)dx\\
&=&\frac{1}{2\pi^2 t^{3/2}(t-s)^{3/2}}\int_{\mathbb{R}^2}
\exp\left(-\frac{t|x_1|^2+|x_2|^2}{t^2}-\frac{(t-s)|x_1-y_1|^2+|x_2-y_2|^2}{(t-s)^2}\right)dx\\
&=&K(s,y)\frac{2\pi s^{3/2}}{4\pi^2 t^{3/2}3(t-s)]^{3/2}}\\
&&\times\int_{\mathbb{R}^2}
\exp\left(\frac{s|y_1|^2+|y_2|^2}{s^2}-\frac{t|x_1|^2+|x_2|^2}{t^2}-\frac{(t-s)|x_1-y_1|^2+|x_2-y_2|^2}{(t-s)^2}\right)dx.
   \eess
Since
   \bess
&&\frac{|y_1|^2}{s}-\frac{|x_1|^2}{t}-\frac{|x_1-y_1|^2}{t-s}\\
&\geq&\frac{|y_1|^2}{s}-\frac{|x_1-y_1|^2+|y_1|^2+2|x_1-y_1||y_1|}{t}-\frac{|x_1-y_1|^2}{t-s}\\
&=&\frac{1}{t}\left(-2|x_1-y_1||y_1|+\frac{t-s}{s}|y_1|^2\right)-\frac{|x_1-y_1|^2}{t}-\frac{|x_1-y_1|^2}{t-s}\\
&\geq&-\frac{s|x_1-y_1|^2}{t(t-s)}-\frac{|x_1-y_1|^2}{t}-\frac{|x_1-y_1|^2}{t-s}\\
&\geq& -\frac{2|x_1-y_1|^2}{t-s}\ \ {\rm for}\ 0<s<t,
  \eess
  and
    \bess
&&\frac{|y_2|^2}{s^2}-\frac{|x_2|^2}{t^2}-\frac{|x_2-y_2|^2}{(t-s)^2}\\
&\geq&\frac{|y_2|^2}{s^2}-\frac{|x_2-y_2|^2+|y_2|^2+2|x_2-y_2||y_2|}{t^2}-\frac{|x_2-y_2|^2}{(t-s)^2}\\
&=&\frac{1}{t^2}\left(-2|x_2-y_2||y_2|+\frac{t^2-s^2}{s^2}|y_2|^2\right)-\frac{|x_2-y_2|^2}{t^2}-\frac{|x_2-y_2|^2}{(t-s)^2}\\
&\geq&-\frac{s^2|x_2-y_2|^2}{t^2(t^2-s^2)}-\frac{|x_2-y_2|^2}{t^2}-\frac{|x_2-y_2|^2}{(t-s)^2}\\
&\geq& -\frac{2|x_2-y_2|^2}{(t-s)^2}\ \ {\rm for}\ 0<s<t,
  \eess
we get for $0<s<t$
   \bess
&&\int_{\mathbb{R}^2}
\exp\left(\frac{s|y_1|^2+|y_2|^2}{s^2}-\frac{t|x_1|^2+|x_2|^2}{t^2}-\frac{(t-s)|x_1-y_1|^2+|x_2-y_2|^2}{(t-s)^2}\right)dx\\
&\geq& \int_{\mathbb{R}^2}
\exp\left(-\frac{2|x_1-y_1|^2}{t-s}-\frac{2|x_2-y_2|^2}{(t-s)^2}\right)dx\\
&=&C_3(t-s)^{3/2}.
  \eess
Substituting the above estimate into (\ref{3.3}) and applying Jensen's inequality, we obtain
   \bess
G(t)&\geq&\frac{C_2}{t^{3/2}}+C_3\int_0^t\left(\frac{s^{3/2}}{t^{3/2}}\right)^pG^p(s)ds \ \ {\rm for}\ t>1.
   \eess
We can rewrite the above inequality as
   \bes
t^{3p/2} G(t)&\geq& C_2t^{3(p-1)/2}+C_3\int_0^ts^{3p/2}G^p(s)ds\lbl{3.4}\\
&=:&g(t).\nm
   \ees
Then for $t>1$, we have
   \bess
&&g(t)\geq C_2t^{3(p-1)/2},\\
&&g'(t)\geq C_3t^{3p/2}G^p(t)\geq C_3t^{3p/2}\left(\frac{1}{t^{3p/2}}g(t)\right)^p=C_3t^{3p(1-p)/2}g^p(t),
  \eess
which implies
   \bess
\frac{C_2^{1-p}}{p-1}t^{-3(p-1)^2/2}\geq \frac{1}{p-1}g^{1-p}(t)\geq C_3\int_t^Ts^{3p(1-p)/2}dx\ \ {\rm for}\
T>t\geq1.
   \eess
If $p\leq1+\frac{2}{3p}$, the right-hand side of the above inequality is
 unbounded as $T\to\infty$, which gives a contradiction in this case. In the
 case $1+\frac{2}{3p}<p<\frac{5}{3}$, we have $3(p-1)^2/2>-1+3p(p-1)/2$, thus
 we get a contradiction by letting $T\to\infty$ and then taking $t\gg1$.

In the case $p=\frac{5}{3}$, we derive from (\ref{3.2}), for $t>1$,
   \bess
u^p(t,x)\geq I_1^p(t,x)\geq\frac{C_1^p}{t^{3p/2}} \exp\left(-\frac{2px_1^2}{t}-\frac{2px_2^2}{t^2}\right).
  \eess
Substituting this estimate into the expression of $I_2$, we obtain, for $t>2$,
   \bess
u(t,x)&\geq& I_2(t,x)\\
&\geq&\int_1^t\int_{\mathbb{R}^2}K(t-s,x-y)\frac{C_1^{p}}{ s^{3p/2}}\exp\left(-\frac{2py_1^2}{s}-\frac{2py_2^2}{s^2}\right)dyds\\
&\geq&\frac{C_4}{t^{3/2}}\exp\left(-\frac{t|x_1|^2+|x_2|^2}{t^2}\right)
\int_1^{t/2}\frac{t^{3/2}}{s^{5/2}(t-s)^{3/2}}ds\\
&&\times\int_{\mathbb{R}^2}\exp\left(\frac{t|x_1|^2+|x_2|^2}{t^2}
-\frac{(t-s)(|x_1|^2+|y_1|^2)+|x_2|^2+|y_2|^2}{(t-s)^2}-\frac{2ps|y_1|^2+2p|y_2|^2}{s^2}
\right)dy\\
&\geq&\frac{C_4}{t^{3/2}}\exp\left(-\frac{t|x_1|^2+|x_2|^2}{t^2}\right)
\int_1^{t/2}\frac{t^{3/2}}{s^{5/2}(t-s)^{3/2}}ds\\
&&\times\int_{\mathbb{R}^2}\exp\left(\frac{t|x_1|^2+|x_2|^2}{t^2}
-\frac{t(|x_1|^2+|y_1|^2)+|x_2|^2+|y_2|^2}{t^2}-\frac{2ps|y_1|^2+2p|y_2|^2}{s^2}
\right)dy\\
&\geq&\frac{C_5}{t^{3/2}}\exp\left(-\frac{t|x_1|^2+|x_2|^2}{t^2}\right)
\int_1^{t/2}\frac{ds}{s}\\
&=&\frac{C_5}{t^{3/2}}\exp\left(-\frac{t|x_1|^2+|x_2|^2}{t^2}\right)\log(t/2).
   \eess
Therefore, for $t>2$, we have
   \bess
G(t)\geq\int_\mathbb{R}K(t,x)\frac{C_5}{t^{3/2}}
\exp\left(-\frac{t|x_1|^2+|x_2|^2}{t^2}\right)\log(t/2)dx
\geq \frac{C_6}{t^{3/2}}\log(t).
   \eess
Using the above estimate, we obtain from (\ref{3.4}), for $t>2$,
   \bess
t^{3p/2}G(t)=\frac{1}{2}t^{3p/2}G(t)+\frac{1}{2}t^{3p/2}G(t)\geq C_7t\log(t)+\frac{C_3}{2}\int_0^ts^{3p/2}G(s)^pds.
  \eess
Denoting the right-hand side of the above inequality by $g(t)$, we have
   \bess
&&g(t)\geq C_7t\log(t)\\
&&g'(t)\geq \frac{C_3}{2}t^{3p/2}G(t)^p\geq\frac{C_3}{2}t^{3p(1-p)/2}g(t)^p,
  \eess
which implies that
   \bess
\frac{3}{2}C_7^{-\frac{2}{3}}[t\log(t)]^{-\frac{2}{3}}\geq\frac{3}{2}g^{-\frac{1}{2}}(t)
\geq\frac{C_4}{2}\int_t^Ts^{-5/3}ds=\frac{C_4}{2}[t^{-\frac{2}{3}}-T^{-\frac{2}{3}}].
  \eess
Letting $T\to\infty$ and $t\gg1$, we get a contradiction. The proof of (ii) is complete.

Using (\ref{3.1}), we have
  \bess
\inf_{x\in\mathbb{R}^2}u(t,x)\geq\mu+\int_0^t\left(\inf_{x\in\mathbb{R}^2}u(s,x)\right)^pds.
   \eess
Then it is easy to see that the solution of (\ref{2.1}) will blow up in finite time.
And thus we complete
the proof. $\Box$

\medskip

\noindent {\bf Acknowledgment} The first author was supported in part
by NSFC of China grants 11771123, 11531006.

 \end{document}